\documentclass[11pt]{article}
\usepackage{amsmath,amsfonts,latexsym,amssymb,amsthm}
\usepackage{url}

\begin{document}

\title{Smooth values of some quadratic polynomials}
\author{Filip Najman \medskip \\University of Zagreb, Croatia}
\date{}
\maketitle
\begin{abstract}
In this paper, using a method of Luca and the author, we find all values $x$ such that the quadratic polynomials $x^2+1,$ $x^2+4,$ $x^2+2$ and $x^2-2$ are $200$-smooth and all values $x$ such that the quadratic polynomial $x^2-4$ is $100$-smooth.
\end{abstract}
\section{Introduction.}
For any integer $n$ we let $P(n)$ be the largest prime factor of $n$ with the convention $P(0)=P(\pm 1)=1$. We say that a integer $s$ is $m$-smooth if $P(s)<m$. In 1964, Lehmer \cite{leh} found all positive integer solutions $x$ to the inequality $P(x(x+1))\le 41$. Notice that 
this amounts to finding all odd positive integers $y=2x+1$ such that $P(y^2-1)\le 41$. There are $869$ such solutions. In \cite{Lu}, Luca found all positive integer solutions of the similar looking inequality  $P(x^2+1)<100$. There are $156$ of them. In \cite{san}, Guzm\'an S\'anchez found the largest solutions of $P(x^2+2)<100$ and $P(x^2-2)<100$. Recently, Luca and the author, in \cite{ln} extended Lehmer's results, finding all the solutions to $P(x^2-1)<100$. There are $16167$ solutions of this inequality.

\medskip

In this paper, we find all the positive integer $x$ that are a solution to the inequalities \begin{equation}P(x^2+1)<200, \label{plus1}\end{equation}
extending Luca's results, \begin{equation}P(x^2+4)<200,\label{plus4}\end{equation}\begin{equation}P(x^2-4)<100,\label{minus4}\end{equation} 
\begin{equation}P(x^2+2)<200,\label{plus2}\end{equation} and \begin{equation}P(x^2-2)<200.\label{minus2}\end{equation}
By solving (\ref{plus2}) and (\ref{minus2}) we will extend the results of Guzm\'an S\'anchez (see \cite{san}).

When considering the inequality (\ref{minus4}), we have to consider all the primes in the factorization of $x^2-4$ (except $2$, as we will later see), while when considering (\ref{plus1}) and (\ref{plus4}) we can disregard the primes congruent to $3$ modulo $4$ and when considering (\ref{plus2}) and (\ref{minus2}) we need only consider the prime $2$ and the primes congruent to $1$ or $3$ modulo $8$ for (\ref{plus2}) and the primes congruent to $1$ or $7$ modulo $8$ for (\ref{minus2}).  This is the reason we can only find $100$-smooth values of $x^2-4$, compared to $200$-smooth values of all the other quadratic forms.

\medskip

In \cite{leh}, \cite{Lu} and \cite{ln} the following approach was taken. Assume that $x$ is a positive integer such that 
$P(x^2\pm 1)\le K$ for the appropriate $K$. Then we can write
\begin{equation}
\label{eq:2}
x^2\pm 1=dy^2,
\end{equation}
where $d$ is squarefree, and $P(dy)\le K$. This implies that only a small number $n$ of primes participates in the factorization of $d$. For example, $n=13$ for $P(x^2-1)\leq 41$ (as in \cite{leh}), $n=12$ for $P(x^2+1)<100$ (as in \cite{Lu}) and $n=25$ for $P(x^2-1)<100$ (as in \cite{ln}).
In all the cases, we can write equation \eqref{eq:2} in the form
$$
x^2-dy^2=\mp 1.
$$
Thus, our possible values for $x$ appear as the first coordinate of one of the solutions of at most $2^{n}-1$ Pell equations. For a given Pell equation, 
the sequence $(y_n/y_1)_{n\ge 1}$  forms a Lucas sequence with real roots. The {\it Primitive Divisor Theorem} for Lucas sequences with real roots (see, for example, \cite{car}, or the more general result from \cite{BHV} which applies to all Lucas sequences) says that if $n>6$, then $y_n$ has a prime factor which 
is at least as large as $n-1$. In the mentioned papers it suffices to check the first $42,\ 98$ and $98$ values respectively of the component $x$ of the Pell equations involved and among these one finds all possible solutions of the equations considered.

\medskip

It is easy to see that in these calculations the number of equations becomes huge. An even bigger problem are the coefficients of the Pell equations, as the size of the solutions grows exponentially in respect to $d$. This means that just \emph{writing down} the solution takes exponential time. This is why in \cite{ln} \emph{compact representations} of the solutions to the Pell equations were used. A compact representation of an algebraic number $\beta \in \mathbb Q (\sqrt d)$ is a representation of $\beta$ of the form
\begin{equation}
\beta=\prod_{j=1}^k \left(\frac{\alpha_j}{d_j}\right)^{2^{k-j}}, 
\label{cr}
\end{equation}
where 
$d_j\in
\mathbb Z,\ \alpha_j={(a_j+b_j\sqrt d)}/2 \in \mathbb Q(\sqrt
d),\ a_j,b_j\in \mathbb Z,\ j=1,\ldots, k$, and $k$, $\alpha$ and $d_j$ have $O(\log d)$ digits. A detailed description of compact representations and their use can be found in \cite{jw}. Using compact representation cuts down the space needed from exponential to polynomial in respect to $d$ and the time needed to compute the regulator of the appropriate real quadratic field from exponential to subexponential. Without the use of compact representations, the results of \cite{ln} would be hopelessly unattainable. For comparison, in \cite{leh} 8191 equations were considered, the
largest of them having $d = 304250263527210$ and the largest solution having less than $8600$ digits, while in \cite{ln} $33554431$ equations were considered, the largest of them having $d=2305567963945518424753102147331756070$, and many of the solutions having billions of digits.

\medskip

In this paper, we will use the same strategy to improve on the results of \cite{Lu} concerning the inequality (\ref{plus1}), using compact representations where standard representations were previously used. In examining the inequalities (\ref{plus4}) and (\ref{minus4}) we will consider fundamental units $\eta_d$ of the quadratic fields $\mathbb Q(\sqrt d)$ satisfying the condition that the cube of $\eta_d$ is the fundamental solution of the Pell equation $x^2-dy^2=-1$ for the case (\ref{plus4}) and $x^2-dy^2=1$ for the case (\ref{minus4}). When considering the inequalities (\ref{plus2}) and (\ref{minus2}), we will use the fact that if $\nu$ is the smallest solution of the equation $x^2-dy^2=\pm 2$, then $\nu^2/2=\eta_d$. This means that all the solutions of $x^2-dy^2=\pm 2$ are of the form $\nu^{2n+1}/2^n$. One can show that the solutions form a \emph{Lehmer sequence} with real roots, and by a result of Ward (see \cite{mw}, and \cite{BHV} for a more general result), we again have control over the prime factors of the members of the sequence. Again, it is necessary to use compact representations when representing the solutions of all these equations.

\medskip

Note that $k=\pm 1 ,\pm 2, \pm4$ are the only values of $k$ such that the integer solutions $x$ to the inequality $P(x^2+k)<M$, for some bound $M$, can be determined in this manner. This is because only for these cases do the solutions of the corresponding equations form either a Lucas or Lehmer sequence. Buchmann, Gy\H ory, Mignotte and Tzanakis found all solutions of $P(x^3+1)<31$, which leads to $P(x^2+3)<31$ in  \cite{bgmt}. They did this by a case by case study and considered each of the Pellian equations (there are $16$ of them) separately.

\section{The inequality $P(x^2+1)<200$.}

Let $x$ be a integer such that $x^2+1$ is $200$-smooth. By a classical result of Fermat, a number can be represented as a sum of to squares only if it has no prime divisors congruent to $3$ modulo $4$. This means that $x^2+1$ can be divisible only by $2$ and the $21$ primes congruent to $1$ modulo $4$ up to $200$. We can now write $x^2+1=dy^2$, i.e. $x^2-dy^2=-1$, where
$$d=2^{a_1}\cdot 5^{a_2}\cdots 197^{a_{22}},\ a_i\in \{0,1\}\text{ for }i=1, \ldots 22.$$
This is a negative Pell equation, so $(x,y)=(X_n,Y_n)$, where $X_n+Y_n\sqrt d= (X_1+Y_1\sqrt d)^{2n+1}$, for some positive integer $n$ and  $X_1+Y_1\sqrt d$ being the fundamental solution of the negative Pell equation. We have $2^{22}-1=4194303$ equations that we need to consider. The largest $d$ appearing is $d=940258296925944608662895221235664431210$. Note that this $d$ is more than $400$ times larger than any $d$ appearing in \cite{ln}.

One can easily check that $Y_1$ divides $Y_n$. We define
$u_n=\frac{Y_n}{Y_1}$. The sequence $(u_n)$ is a Lucas sequence of the first kind with real roots $\eta$ and $\zeta$. By a result of Carmichael (see \cite{car}), the Primitive Divisor Theorem, this implies that for every $n>12, u_n$ has a \emph{primitive divisor} $p$, a prime satisfying, among other properties, $p\equiv\pm 1 \pmod n$. This implies that for $n>200$, there exists a prime $p\geq 200$ dividing $u_n$. Thus, only the first $200$ values of $(Y_n)$ can possibly be $200$-smooth. This means that we only need to consider the first $200$ solutions of each negative Pell equation.

We will follow the methods of \cite{ln} very closely, so we will give only an outline of the algorithm (for a detailed explanation of each step see \cite{ln}). The algorithm that finds all the solutions goes through all the possible values $d$ and for each, does the following: First it computes the regulator of the quadratic field $\mathbb Q(\sqrt d)$ using Buchmann's subexponential algorithm (see \cite{bc}). The results of this algorithm are dependent on the Generalized Riemann hypothesis. Next, using the computed regulator, we construct a compact representation of the fundamental solution $X_1+Y_1\sqrt d$ of the negative Pell equation $x^2-dy^2=-1$, using the methods from \cite{mm}. Now, using the algorithm for modular arithmetic described in \cite{fn}, we check whether $Y_1$ is $200$-smooth. If it is, then we check which of the $Y_n,\ 2\le n\le 200$ is $200$-smooth. Each $200$-smooth value of $Y_n$ gives us a solution $X_n$ of our problem. If $Y_1$ is not $200$-smooth, then no $Y_n$ will be. There is one last check that needs to be done: compute all the convergents $\frac{p_n}{q_n}$ of the continued fraction expansion of $\sqrt d$ such that $q_n<z$, where $z$ is the $200$-smooth part of $Y_1$, and check whether there exists a convergent $\frac{p_n}{q_n}$ such that $p_n^2-q_n^2=-1$. This test removes the dependence of our results on the Generalized Riemann hypothesis. In all tested cases it failed, as if a convergent $\frac{p_n}{q_n}$ satisfying the test were to be found, it would imply that the Generalized Riemann hypothesis is false.

We obtain the following results:
\newtheorem{tm}{Theorem}
\begin{tm} 

\begin{itemize}

\item[a)] The largest three solutions of the equation $P(x^2+1)<200$ are
$x=69971515635443
$, $120563046313$ and\\ $104279454193$.
\item[b)] The largest solution of $P(x^4+1)<200$ is $x=10$.
\item[c)] The largest solution of $P(x^6+1)<200$ is $x=8$.
\item[d)] The largest $n$ such that $P(x^{2n}+1)<200$ has a solution is $n=9$, the solution being $x=2$.
\item[e)] The inequality $P(x^2+1)<200$ has $811$ solutions.
\item[f)] The greatest power $n$ of the fundamental solution of the negative Pell equation $(X_1+Y_1\sqrt d)^{n}$ which leads to a solution of our problem is $n=9$ for $d=5$. The case $d=5$ also gives us most solutions, namely $4$ of them.

 \end{itemize}
\end{tm}
\emph{Proof:}\\
Part b) is proved by finding the largest square of all the $x$, c) by finding the largest a cube, etc. \qed \\

\section{The inequality $P(x^2+4)<200$.}

Obviously, $P(x^2+1)<200$ iff $P((2x)^2+4)<200$. Thus, we have allready obtained all the even solutions to the inequality (\ref{plus4}). It remains to find the odd solutions. Let $\eta_d=\frac{u+v\sqrt d}{2}$ be the fundamental unit of the quadratic field $\mathbb Q(\sqrt d)$ and $x_1+y_1\sqrt d$ the fundamental solution of the Pell equation $x^2-dy^2=1$. Then $x_1+y_1\sqrt d=\eta_d^n$, where $n=1,2,3$ or $6$, and the exact value can be found by examining $u$ and $v$ modulo $8$, from the following table: 
\begin{center}
\begin{tabular}{|c|c|c|c|}
\hline
$d$& $v$& $u$&$n$\\
\hline
$d\equiv 1 \pmod 4$ & $v\equiv 0 \pmod 4$ & $-$ & $1$\\
 & $v\equiv 2 \pmod 4$ & $-$ & $2$\\
\hline
$d\equiv 5 \pmod {16}$ & $v\equiv 1 \pmod 2$ & $u\equiv \pm 3 v\pmod 8$ & $3$\\
$$ & $$ & $u\equiv \pm  v\pmod 8$ & $6$\\
\hline
$d\equiv 13 \pmod {16}$ & $v\equiv 1 \pmod 2$ & $u\equiv \pm  v\pmod 8$ & $3$\\
$$ & $$ & $u\equiv \pm  3v\pmod 8$ & $6$\\
\hline
$d\equiv 2 \pmod 4$ & $v\equiv 0 \pmod 2$ & $-$ & $1$\\
 & $v\equiv 1 \pmod 2$ & $-$ & $2$\\
\hline
$d\equiv 3 \pmod 4$ & $-$ & $-$ & $1$\\
\hline
\end{tabular}
\bigskip\\
Table 1\\
\end{center}
It is easy to see that $x^2-dy^2=-4$ will have a solution iff $n=6$ in Table 1, and then the solutions will be of the form $$X_m+Y_m\sqrt d=2\eta_d^k,$$ where $k\equiv \pm 1 \pmod 6$. Also, one can see that $n=6$ only if $d\equiv 5 \pmod 8$. This also means that we do not have to consider $2$ in the factorization of $d$, leaving only the $21$ primes congruent to $1$ modulo $4$. 

Using the same algorithm from Section 2, with the appropriate minor changes, we obtain the following results. 
\newtheorem{tm2}[tm]{Theorem}
\begin{tm2} 

\begin{itemize}

\item[a)] The largest three odd solutions of the equation $P(x^2+4)<200$ are
$x= 191686681859$, $112899039159$ and $28608252345 $.
\item[b)] The largest odd solution of $P(x^4+4)<200$ is $x=923$.
\item[c)] There are no odd solutions to $P(x^{2n}+4)<200$ for $n\geq 3$.
\item[d)] The inequality $P(x^2+4)<200$ has $344$ odd solutions.
\item[f)] The greatest power $n$ of the smallest solution of our problem $(X_1+Y_1\sqrt d)^{n}$ which leads to a solution  is $n=19$ for $d=5$. The case $d=5$ also gives the most solutions, namely $5$ of them. \end{itemize}
\end{tm2}

\section{The inequality $P(x^2-4)<100$.}

Obviously, $P(x^2-1)<100$ iff $P((2x)^2-4)<100$. Thus, we have already obtained in \cite{ln} all the even solutions to the inequality (\ref{plus4}). It remains to find the odd solutions. Again, let $\eta_d=\frac{u+v\sqrt d}{2}$ be the fundamental unit of the quadratic field $\mathbb Q(\sqrt d)$ and $x_1+y_1\sqrt d$ the fundamental solution of the Pell equation $x^2-dy^2=1$, and $x_1+y_1\sqrt d=\eta_d^n$, where $n$ can be found in Table 1. It is easy to see that $x^2-dy^2=4$ will have a solution iff $n=3$ in Table 1, and then the solutions will be of the form $$X_m+Y_m\sqrt d=2\eta_d^k,$$ where $k\not\equiv 0 \pmod 3$. Again, one can see that $n=3$ only if $d\equiv 5 \pmod 8$, so once again we can disregard the prime $2$ in the factorization of $d$.lready obtained in \cite{ln} all the even solutions to the inequality (\ref{plus4}). It remains to find the odd solutions. Again, let $\eta_d=\frac{u+v\sqrt d}{2}$ be the fundamental unit of the quadratic field $\mathbb Q(\sqrt d)$ and $x_1+y_1\sqrt d$ the fundamental solution of the Pell equation $x^2-dy^2=1$, and $x_1+y_1\sqrt d=\eta_d^n$, where $n$ can be found in Table 1. It is easy to see that $x^2-dy^2=4$ will have a solution iff $n=3$ in Table 1, and then the solutions will be of the form $$X_m+Y_m\sqrt d=2\eta_d^k,$$ where $k\not\equiv 0 \pmod 3$. Again, one can see that $n=3$ only if $d\equiv 5 \pmod 8$, so once again we can disregard the prime $2$ in the factorization of $d$.

Using the same algorithm from Section 2, with the appropriate minor changes, we obtain the following results. 
\newtheorem{tm3}[tm]{Theorem}
\begin{tm3} 

\begin{itemize}

\item[a)] The largest three odd solutions of the equation $P(x^2-4)<100$ are
$x=407479035814853$, $335682488669673$ and\\ $250734674482437$.
\item[b)] The largest odd solution of $P(x^4-4)<100$ is $x=59$.
\item[c)] The largest odd solution of $P(x^6-4)<100$ is $x=7$.
\item[d)] The largest $n$ such that $P(x^{2n}-4)<100$ has a solution is $n=7$, the solution being $x=3$.
\item[e)] The inequality $P(x^2-4)<100$ has $2846$ odd solutions.
\item[f)] The greatest power $n$ of the smallest solution of our problem $(X_1+Y_1\sqrt d)^{n}$ which leads to a solution of our problem is $n=10$ for $d=5$. The case $d=5$ also gives the most solutions, namely $7$ of them. \end{itemize}
\end{tm3}

\section{The inequality $P(x^2+2)<200$.}
\label{ch5}
Consider the equation $x^2-dy^2=\pm2$, where $d$ is square free, and let $x_1+y_1\sqrt d=\nu$ be the least solution. Obiviously, $d\equiv 2,3 \pmod 4$ has to hold. One can see that, using the notation from Table 1, this implies $n=1$ or $2$. From a result of Perron (see \cite[p. 126-129]{per}), if $d\neq 2$, at most one of the equations $x^2-dy^2=-1$, $x^2-dy^2=2$ and $x^2-dy^2=-2$ is solvable. This implies that $n=2$ is impossible. By a elementary argument (see for example \cite[p. 420]{jw}), one can also show that $\nu^2/2=\eta_d$ holds. This can be reformulated in terms of the infrastructure of a real quadratic field by saying the ideal with norm $\pm 2$ appears half way through the cycle of the principal class, or in terms of continued fractions by saying that the $p_{l/2}^2-dq_{l/2}^2=\pm 2$, where $l$ is the length of the continued fraction expansion of $\sqrt d$. Thus all the solutions of $x^2-dy^2=\pm2$ are of the form $$x_k+y_k\sqrt d = \nu\eta_d^k=\frac{\nu^{2k+1}}{2^k},\ \text{ for } k\leq 0.$$  
Now let $$\alpha=\frac{x_1+y_1\sqrt{d}}{\sqrt 2}\text{ and } \beta=\frac{x_1-y_1\sqrt{d}}{\sqrt{2}}.$$
Then $$y_k=\frac{\alpha^k-\beta^k}{\sqrt{2d}},$$ and
$$u_k=
\begin{cases}
\frac{y_k}{y_1} \text{ if $k$ is odd}\\
\frac{y_k}{y_2} \text{ if $k$ is even}\\
\end{cases}
$$
is a Lehmer sequence. Note that only the odd members of this sequence yield a solution to our equation. Thus, by the Primitive divisor theorem for Lehmer sequences (see \cite{mw} and \cite{BHV}), only the first $100$ odd members of the sequence $(u_k)_{k\geq 0}$ can possibly be $200$-smooth.

\medskip

To narrow down our search we will use the following results of Yokoi:

\newtheorem{tm4}[tm]{Theorem} 
\begin{tm4}{\bf (\cite[Theorems 1 and 2]{yo})}\\
Let $t+u\sqrt d$, where $d$ is a squarefree positive integer congruent to $2$ or $3$ modulo $4$, be the fundamental unit of the real quadratic field $\mathbb Q(\sqrt d)$.\\
a) The Diophantine equation $x^2-dy^2=2$ is solvable if and only if $t\equiv 1 \pmod d$.\\
b) The Diophantine equation $x^2-dy^2=-2$ is solvable if and only if $t\equiv -1 \pmod d$.\\    
\label{te4}
\end{tm4}
\newtheorem{tm5}[tm]{Proposition} 
\begin{tm5}{\bf (\cite[Proposition 2]{yo})}\\
a) If the Diophantine equation $x^2-dy^2=2$ is solvable then $p\equiv \pm 1 \pmod 8$ for any odd prime factor $p$ of $d$.\\
b) If the Diophantine equation $x^2-dy^2=-2$ is solvable then $p\equiv 1\text{ or } 3 \pmod 8$ for any odd prime factor $p$ of $d$.\\   
\label{te5}
\end{tm5}

Using Theorem \ref{te4} and Proposition \ref{te5} to narrow down our search, we procced in the same manner as for (\ref{plus1}), (\ref{plus4}) and (\ref{minus4}). We obtain the following results:

\newtheorem{tm6}[tm]{Theorem}
\begin{tm6} 

\begin{itemize}

\item[a)] The largest three odd solutions of the equation $P(x^2+2)<200$ are
$x=9575480365630$, $14629598023$ and $8850900308$.
\item[b)] The largest odd solution of $P(x^4+2)<200$ is $x=171$.
\item[c)] The largest solution of $P(x^6+2)<200$ is $x=3$.
\item[d)] The largest $n$ such that $P(x^{2n}+2)<200$ has a solution is $n=5$, the solution being $x=2$.
\item[e)] The inequality $P(x^2+2)<200$ has $914$ solutions, of which $516$ are odd and $398$ are even.
\item[f)] The greatest power $n$ such that $$\frac{(X_1+Y_1\sqrt d)^{2n+1}}{2^n}$$ leads to a solution of our problem is $n=8$ for $d=2$. The case $d=3$ gives the most solutions, namely $6$ of them. \end{itemize}
\end{tm6}

\section{The inequality $P(x^2-2)<200$.}
We search for the solutions of this inequality as explained in the previous section, again using Theorem \ref{te4} and Proposition \ref{te5} to narrow down our search. We obtain the following results:   

\newtheorem{tm7}[tm]{Theorem}
\begin{tm7} 

\begin{itemize}

\item[a)] The largest three  solutions of the equation $P(x^2-2)<200$ are
$x=324850200677887$, $1600947755823$ and $494400410248$.
\item[b)] The largest solution of $P(x^4-2)<200$ is $x=47$.
\item[c)] The largest odd solution of $P(x^6-2)<200$ is $x=10$.
\item[d)] The largest $n$ such that $P(x^{2n}-2)<200$ has a solution is $n=5$, the solution being $x=2$.
\item[e)] The inequality $P(x^2-2)<200$ has $537$ solutions, of which $313$ are even and $224$ are odd.
\item[f)] The greatest power $n$ such that  $$\frac{(X_1+Y_1\sqrt d)^{2n+1}}{2^n}$$ leads to a solution of our problem is $n=7$ for $d=2$. The case $d=2$ also gives the most solutions, namely $5$ of them. \end{itemize}
\end{tm7}

Problem number $4$ on the list \cite{et} of open problems concerning Diophantine equations is to find all the solutions $(x,y,p)$ to 
$$x^2-2=y^p,$$
where $p$ is a odd prime. By searching through our solutions we prove the following proposition.
\newtheorem{tm8}[tm]{Proposition} 
\begin{tm8}
If the Diophantine equation $x^2-2=y^p$, where $p$ is an odd prime, then $P(y)>200$.\\
 
\end{tm8}

\textbf{Remark.}
The tables produced by our computations can be found on the web page \url{http://web.math.hr/~fnajman}. 
\bigskip

\textbf{Acknowledgements.}
Many thanks go to Florian Luca for pointing out that this method will work for $x^2\pm 2$. I am also grateful to Andrej Dujella and the referee for many helpful comments.

\bigskip

\noindent\textit{Department of Mathematics, University of Zagreb, Bijeni\v cka cesta 30, 10000 Zagreb, Croatia\\
fnajman@math.hr}
\end{document}